\documentclass[12pt]{article}
 \usepackage{amsfonts,amssymb,graphicx}
\usepackage[pdfpagemode=None,colorlinks=true,linkcolor=blue,citecolor=blue]{hyperref}

\newtheorem{prop}{Proposition}
\newtheorem{rema}{Remark}

\newtheorem{lemm}{Lemma}

\newtheorem{coro}{Corollary}
\newtheorem{exem}{Example}

\newcommand{\R}[1][]{\ensuremath{{\mathbb{R}^{#1}} }}
\renewcommand{\S}[1][]{\ensuremath{{\mathbb{S}^{#1}} }}

\newcommand{\<}{\langle}
\renewcommand{\>}{\rangle}
\newcommand{\ga}{\gamma}
\newcommand{\pa}{\partial}

\newcommand{\eps}{\epsilon}
\newcommand{\te}{\theta}

\date{}

\title{The Blaschke-Lebesgue problem for constant width bodies of revolution}
\author{Henri Anciaux\footnote{The first author is supported by SFI (Research Frontiers Program)} , Nikos Georgiou}
\begin{document}
\maketitle

\centerline{\textbf {Abstract}}

\smallskip

{\small
We prove that among
 all constant width
 bodies of revolution, the minimum of the ratio of the volume to the cubed width
is attained by
the constant width body
obtained by rotation of the Reuleaux triangle about an axis of symmetry.}
\medskip

\centerline{\small \em 2000 MSC: 52A15
\em }


\section*{Introduction}
The width of a convex body $B$ in $n$-dimensional Euclidean
space in the direction $\vec{u}$ is the distance between the two
supporting planes of $B$ which are orthogonal to $\vec{u}.$ When
this distance is independent of $\vec{u}$, $B$ is said to have \em
constant width. \em

The ratio ${\cal I}(B)$ of the volume of a constant width body to 
the volume of the ball of the same width
 is homothety invariant, as is the isoperimetric ratio. Moreover the maximum of
 ${\cal I }(B)$ is
 attained by  round spheres, just as the minimum of the isoperimetric ratio.
 However, while the latter is not bounded from above, the infimum of ${\cal I}$ is strictly
 positive, since
compactness properties of the space of convex sets ensures the
existence of a minimizer.  It is known
since the works of Blaschke and Lebesgue that the \em Reuleaux triangle, \em
obtained by taking the intersection of three discs centered at the
vertices of an equilateral triangle, minimizes ${\cal I}$ in
dimension $n=2.$ The determination of the minimizer of ${\cal I}$ in any dimension is
the  \em Blaschke-Lebesgue problem. \em Recently several simpler solutions of the problem in dimension $2$
have been given (cf \cite{Ba},\cite{Ha}), however the Blaschke-Lebesgue
problem in dimension $n=3$ appears to be very difficult to solve
and remains open.

\medskip

In this paper we prove:

\bigskip

\noindent \textbf{Main Theorem:} \em Amongst all constant width
 bodies of revolution in Euclidean $3$-space, the minimum of the ratio the volume to the cubed width
is attained by
the constant width body $B_{Reul}$
obtained by rotation of the Reuleaux triangle about an axis of symmetry. \em

\bigskip

This result has been proved in \cite{CCG} by a geometric argument
 which relies on a direct comparison
between the volume of an arbitrary convex body of revolution of constant width, and
 the rotated Reuleaux triangle of the same width. Our proof is more
analytical in nature as it uses calculus of variations. Moreover,
some of the observations we make do not depend on the assumption of
rotational  symmetry.

It is known that there exists a convex body $B_{Meis},$ called \em Meissner's
tetrahedron, \em satisfying ${\cal I}(B_{Meis}) \simeq 0.8019$ (cf \cite{CG},\cite{GK}, \cite{Ba}),
 while  ${\cal I}(B_{Reul}) \simeq  0.8584$. Thus the main theorem implies the following:

\bigskip

\noindent \textbf{Corollary:} \em The solution of the
Blaschke-Lebesgue problem is not a body of revolution. \em

\bigskip

As in \cite{Ba} and \cite{Ha}, our proof is based on the analysis
of the support function $s$ which characterizes a convex body of
constant width $2w.$ A crucial point is the following
observation, stated in \cite{GK}: flowing the boundary of a convex
body along its inward unit normal vector field preserves the
constant width condition, as long as the evolving surface remains
convex.
 Moreover, the ratio ${\cal I}$ decreases along the flow, so
the minimizer of ${\cal I}$ must occur at the latest time such
that convexity holds, and therefore must be singular. This issue
is easily controlled by  introducing the function $h=s-w$, which
 is invariant along the normal flow, while the width
$2w$ decreases linearly. Thus, there
exists a positive number $w_0(h)$ such that for any $w \geq
w_0(h),$  the function $s=h+w$ is the support function  of some convex body of
constant width $2w$.
 Hence, we can restrict the minimization process to the class of support functions of the form $s=h+w_0(h),$
while all the necessary information is carried by the function
$h.$ The assumption of rotational symmetry made in the present
 article, since it reduces all the involved calculus to one variable, simplifies considerably the exposition;
however, all these facts hold for arbitrary constant width bodies in Euclidean $3$-space.

The next step in the proof of our main theorem
consists of
using the second order condition of minimization (i.e. stability) to prove that the map $|h''+h|$ must be constant. It follows
 that the value of  ${\cal I}(h)$ is completely determined by the set of the
discontinuities of $h''+h$. The rotated Reuleaux triangle $B_{Reul}$ corresponds to the case of
$h''+h$ having the least possible number of discontinuities.
 We then show  that, unless the number of discontinuities of $h''+h$ is minimal, one can always reduce the ratio
${\cal I}$, 
which completes the proof.

 As a final comment, we point out that
one can
 prove the fact that the
Reuleaux triangle minimizes ${\cal I}$ among the constant width
bodies of the plane by a slight modification of our argument.

\section{Preliminaries: constant width bodies of revolution}
Let $B$ be a convex body of revolution in  $\R^3$,
i.e. it is invariant under rotation around some axis. We may assume without
loss of generality that this axis is vertical. Therefore the
boundary of $B$
 can be parametrized by
 $$ \begin{array}{lccc}  X :
 & [a,b] \times \S^1
 &\to&  \R^3 \\
&  (\phi,\te) & \mapsto &  (x(\phi)\cos \te, x(\phi) \sin \te, y(\phi)),
   \end{array}$$
where  $\ga(\phi)=(x(\phi),y(\phi))$ is a parametrized curve
 such that $x(\phi) \geq 0$ and $x(a)=x(b)=0.$
It is known (see \cite{Ho}, \cite{Ba}) that if $B$ has constant width it must be
strictly convex. It follows that the generating curve $\ga$ is also
 strictly convex,
which allows us to reparametrize it   by the angle $t$ made by its unit outward normal vector
$\vec{n}(t)=(\cos t, \sin t)$ with the horizontal plane:
$\ga(t)=(x(t),y(t))$, with $ t \in [-\pi/2,\pi/2].$
This parametrization holds even if the curve $\ga$
 is not regular.

\medskip

Next we  express the constant width assumption of $B$ in
terms of the curve $\ga$. For this purpose it is convenient
 to consider the union of  $\ga$ with its image under reflection through the vertical axis, which gives a
strictly convex, closed, planar curve. This closed curve is parametrized by
$$\ga(t)=(x(t),y(t)):=(-x(\pi-t),y(\pi - t)), \forall t \in [-\pi,-\pi/2] \cup [\pi/2,\pi],$$
and it is then possible to parametrize $\pa B$ by
$$ \begin{array}{lccc}  X :
 & (\S^1 \times \S^1) / \sim
 &\to&  \R^3 \\
& (t,\te) & \mapsto &  (x(t)\cos \te, x(t) \sin \te, y(t)),
   \end{array}$$
where $\sim$ denotes the equivalence relation defined on the torus
$\S^1 \times \S^1$ by $(t,\te) \sim (\pi - t , \te+\pi).$ In
particular, the antipodal point of $(t,\te)$ is $(-t,\te+\pi) \sim
(t+\pi,\te).$

Next  the support function of  $\pa B$ at the point $X(t,\te)$ is
defined to be
$$s_X(t,\te) :=\< X(t,\te),\vec{N}(t,\te) \>,$$
where $\vec{N}$ is the unit outward vector of $\pa B$ at the point $X(t,\te).$
The width of  $B$ in the direction $\pm \vec{N}$ is equal to the sum of the
 support function evaluated at the two antipodal points corresponding to this direction:
$$ 2w(t,\te)=s_X(t,\te)+s_X(-t,\te+\pi)=s_X(t,\te)+s_X(t+\pi,\te).$$
Similarly, the support function of
the curve $\ga$ is
$s_\ga(t):=\< \ga(t),\vec{n}(t)\>=\< \ga(t), (\cos t,\sin t) \>.$
An easy computation, using the fact that $\vec{N}=(\cos t \cos \te, \cos t \sin \te, \sin t),$ gives:
$$ 2w(t,\te)=\<X(t,\te),\vec{N}(t,\te) \>+\<X(t+\pi,\te),\vec{N}(t+\pi,\te) \> $$
$$=\<(x(t),y(t)),(\cos t,\sin t )\>+\<(x(t+\pi),y(t+\pi)),(\cos (t+\pi),\sin (t+\pi )\>$$
$$=s_\ga(t)+s_\ga(t+\pi).$$
The last expression is nothing but the
 width of the curve $\ga$ in the direction $(\cos t, \sin t),$ so we have proved:
\begin{lemm}
$B$ has constant width if and only if $\ga$ has constant width.
\end{lemm}

From now on we focus on the curve $\ga$ and use complex notation
in the Euclidean plane $\{ (x,y) \simeq x+iy \}.$ The curve $\ga$
 can be reconstructed from its support function~$s_\ga$:
 $$\ga(t) = s_\ga(t) e^{it} +  s_\ga'(t) i e^{it}.$$
Differentiating this expression yields
$\ga'(t)=(s_\ga''+s_\ga)ie^{it}$ and thus the curve is regular if
and only if $s_\ga''+s_\ga >0.$ Moreover, it
 changes its convexity with the sign of $s_\ga''+s_\ga.$ As the
 curve may not be regular everywhere, but must remain strictly convex, we are left with the condition $s_\ga''+s_\ga \geq 0$.
 This quantity is nothing but the radius of
 curvature of $\ga.$

\medskip

Next set $w:=\frac{1}{2\pi}\int_{\S^1} s_\ga(t) dt$ and $h:=s_\ga-w.$ Hence,
the curve $\ga$ has constant width if and only if
$$h(t)+h(t+\pi)=0, \leqno{(1)}$$
and in this case the width is exactly $2w.$ On the other hand, the
symmetry of the curve with respect to the vertical axis, i.e.
$x(\pi-t)+iy(\pi-t)=-x(t)+iy(t)$ implies that
$s_\ga(t)-s_\ga(\pi-t)=0$ and thus
$$h(t)-h(\pi-t)=0. \leqno{(2)}$$
By Equations $(1)$ and $(2)$  it is enough to
define $h$ on the interval $[0,\pi/2]$ and to extend it to $\S^1$ by the symmetries $(1)$ and $(2).$
Furthermore we have $h(0)=0$ and $h'(\pi/2)=0$. On the other hand, it is proven in \cite{Ho} that the support
function of a constant width body is $C^{1,1}$  so we conclude that
 the functional space of the $h$ corresponding to  constant width curves with axial symmetry (and thus to constant width bodies
of revolution) is
$$E:= \{  h \in C^{1,1}([0,\pi/2]), h(0)=0, h'(\pi/2)=0 \}.$$
Finally, a given pair $(h,w)$ corresponds to a curve $\ga$ if the
support function $s_\ga=h+w$ satisfies the condition
$s_\ga''+s_\ga =h''+h+ w\geq 0.$ By the Rademacher theorem, the
fact $h \in C^{1,1}$ implies that $h''$ exists a.e. and belongs to
$L^{\infty}(0,\pi/2)$.
 Thus
we must have $w \geq -(h''+h),$ almost everywhere in $\S^1.$ As
$h$ is odd, it is equivalent to require that $w \geq  h''+h,$ a.e. in
$\S^1.$ Hence, for any given $h \in E,$ we define
 $$ w_0(h) := || h(t)+h''(t)||_{L^\infty(0,\pi/2)}.$$

Summing up, we have proven:

\begin{prop} \label{p1} There is a one-to-one correspondence between the convex bodies $B$ of  revolution which have constant width $2w$,
and the set of pairs $(h,w), h \in E, w \geq w_0(h),$ where
$$ E:= \{  h \in C^{1,1}([0,\pi/2]), h(0)=0, h'(\pi/2)=0 \},$$
and
$$  w_0(h) := || h(t)+h''(t)||_{L^\infty(0,\pi/2)}.$$
\end{prop}

\begin{exem}
If $h=c \sin t,$ where $c$ is some real constant,
 the curve $\ga$ is a circle centered in the vertical axis and thus the corresponding body is a ball.
\end{exem}

\section{The Blaschke-Lebesgue problem} \label{s2}
We  now compute the volume of a constant width body of revolution $B$ in terms of $h$ and $w.$
We start by calculating the first derivatives of the immersion:
 $$ X_t=(x' \cos \te ,x' \sin \te, y'), \quad \quad X_\te=(-x \sin \te, x \cos \te,0).$$
As $\det(X,X_\te,X_t)>0,$
the volume of $B$ is
$${\cal V}(B)=\frac{1}{3} \int_0^{2\pi}  \int_{-\pi/2}^{\pi/2} \det(X,X_\te,X_t)   dt d \te$$
$$ = \frac{2\pi}{3}   \int_{-\pi/2}^{\pi/2} x(y'x-yx')   dt.  $$
The integrand turns out to be a polynomial in $w$:
$$x(y'x-yx')=((w+h) \cos t   - h' \sin t)(w+h+h'')(w+h)$$
$$= (w \cos t +h \cos t   - h' \sin t)(w^2+(2h+h'')w+h(h''+h)) $$
$$=  \cos t w^3+  \left(h \cos t - h' \sin t + (2h+h'')\cos t \right) w^2$$
$$+\Big( 2h''h \cos t -h''h' \sin t + 3 h^2 \cos t - 2 hh' \sin t  \Big)w  + (h \cos t - h' \sin t)(h''+h)h .  $$
Next, as $h$ is odd, the functions $h \cos t$, $h'
\sin t$ and $(h+h'') \cos t$ are odd as well, thus the integrals of both
the term in $w^2$, and the constant term, vanish and we are left to compute of the $w$ term.
Using the fact that
 $$ \int_{0}^{\pi/2} \left( 2h'' h' \sin t + (h')^2 \cos t \right)dt=\int_{0}^{\pi/2} \big( (h')^2
 \sin t \big)' dt ,$$
 and
 $$ \int_{0}^{\pi/2} \left(h'' h \cos t -h'h \sin t + (h')^2 \cos t \right) = \int_{0}^{\pi/2} \big( hh' \cos t \big)' dt,$$
 vanish by the boundary conditions, we get
 $$ \int_{-\pi/2}^{\pi/2}(2h''h \cos t -h''h' \sin t + 3 h^2 \cos t - 2 hh' \sin t)dt$$
 $$  =2\int_{0}^{\pi/2}(2h''h \cos t -h''h' \sin t + 3 h^2 \cos t - 2 hh' \sin t)dt$$
 $$ =2\int_{0}^{\pi/2} (3h^2 \cos t -\frac{3}{2}(h')^2 \cos t )dt.$$
Thus we obtain
$${\cal V}(B)  =4 \pi \left( \frac{w^3}{3}  +w   \int_{0}^{\pi/2}(h^2 - \frac{(h')^2}{2})\cos t dt \right).$$

\begin{prop}
Let $(h,w)$ with $h \in E$ and $w \geq w_0(h)$ and let $B$ be the corresponding constant width body of revolution
 (see Proposition \ref{p1}). Then
$${\cal V}( B)= 4\pi \Big(\frac{w^3}{3} + w{\cal F}(h) \Big)$$
 where the functional ${\cal F}$ is defined to be
$${\cal F}(h):=\int_0^{ \pi/2} (h^2- \frac{1}{2}(h')^2)\cos t dt.$$
\end{prop}

\begin{rema}
We could also have computed the area of $\pa B$ and used the Blaschke formula (see \cite{GK}) for bodies of constant width:
${\cal V}(B)={\cal A}(\pa B)w - \frac{8 \pi }{3}w^3. $
\end{rema}
It is then easy to express the ratio ${\cal I}(B)$ in terms of $h$
and $w$:
$${\cal I}( B):=\frac{{\cal V}(B)}{4 \pi w^3/3}= {1}  + \frac{3{\cal F}(h)}{w^2}.$$
The next proposition, which may be seen as a weighted version of
the classical Wirtinger inequality, shows that the last term in
the expression above is negative:

\begin{prop}[weighted Wirtinger inequality] \label{Wirtinger}
Let $h \in E.$ Then the following inequality holds,
$${\cal F}(h)=\int_0^{ \pi/2} (h^2-\frac{1}{2}(h')^2)\cos t dt \leq 0,$$
and the equality is attained if and only if $h=c \sin t$ for some real constant $c.$
\end{prop}

\noindent \textit{Proof.}  Introduce the function $g:=h \cos t - h' \sin t.$ It is easy to
check that $g'= -\sin t(h+h'').$
 The boundary conditions $h(0)=0$ and $h'(\pi/2)=0$ imply that
$g(0)=g(\pi/2)=0,$ and, moreover,
$$\lim_{\eps \to 0 } \frac{g(t)}{\sin t}= \lim_{\eps \to 0 } (h(t)
\cot t - h'(t))=h'(0)-h'(0)=0.$$ Hence
$${\cal F}(h) = \lim_{\eps \to 0 }
\int_\eps^{\pi/2} (h+h'')(h \cos t - h'\sin t)dt$$
$$ =\lim_{\eps \to 0 } \int_\eps^{\pi/2} \frac{gg'}{\sin t}dt
=\lim_{\eps \to 0 }\left( \int_\eps^{\pi/2} \frac{g^2(t)}{2}\left(
\frac{1}{\sin t} \right)'dt + \left[\frac{g^2(t)}{2 \sin t}
\right]_\eps^{\pi/2} \right)$$
$$=\lim_{\eps \to 0 } \left( \frac{-1}{2}\int_\eps^{\pi/2} \frac{g^2(t) \cos t }{\sin^2
t}-\frac{g^2(\eps)}{2 \sin \eps}\right)=
\lim_{\eps \to 0 } \frac{-1}{2}\int_\eps^{\pi/2} \frac{g^2(t) \cos t }{\sin^2 t}\leq 0.$$
The last inequality shows that if ${\cal F}(h)$ vanishes, so does $g.$
We thus have $h \cos t = h' \sin t.$ It is easy to check that the only solutions
of this linear equation with initial conditions $h(0)=0$ and $h'(\pi/2)=0$
are $h= c \sin t,$ where $c$ is some real constant.

\bigskip

From this lemma we recover the fact that the ratio ${\cal I}$
achieves its maximum for $h=c \sin t,$ which corresponds to $B$
being a ball. Moreover, it follows that  the ratio $ {\cal I}(B)$
is increasing with respect to $w.$ So by Proposition \ref{p1}  we
get:
\begin{coro}
Let $(h,w)$ be  a minimizer of ${\cal I}(h,w).$ Then $w=w_0(h).$
\end{coro}

We end this section using again the weighted Wirtinger inequality to prove that $|h+h''|$ must be constant.
\begin{prop} \label{const}
Let $(h,w_0(h))$ be a minimizer of ${\cal I}(h,w).$ Then the
quantity $|h''+h|$ is constant almost everywhere in $[0,\pi/2]$.
\end{prop}

\noindent \textit{Proof.} We proceed by contradiction assuming
that there is a non-empty interval $[a,b]$ included in $[0,\pi/2]$
such that $|h(t)+h''(t)| < w_0(h)$ a. e. in $[a,b]$.
  Consider a map $V$ of $E$ whose support
is contained in $[a,b]$ and which is not of the form $V=c \sin t,$
and  define the deformation $h^\eps:=h+\eps V $ of $h.$ For
small $\eps,$
$$w_0(h^\eps)=|| (h^\eps)'' + h^\eps||_{L^{\infty}(0,\pi/2)} =||h+h''||_{L^{\infty}(0,\pi/2)}=w_0(h),$$
hence
$$\frac{{\cal F}(h^\eps)}{w^2_0(h^\eps)}  =\frac{{\cal F}(h)}{w^2_0(h)} + \eps  \frac{\delta {\cal F}(h,V)}{w^2_0(h)}
+ \frac{\eps^2}{2} \frac{ \delta^2 {\cal F}(h,V)}{w^2_0(h)} +
o(\eps^2).$$ As $h$ is a minimizer of ${\cal I}$, and thus, of
${\cal F}(h)/w^2_0(h)$, we must have both $\delta {\cal F}(h,V)
=0$ and $\delta^2 {\cal F}(h,V) \geq 0.$ On the other hand the
functional ${\cal F}$ is quadratic, so that $\delta^2 {\cal
F}(h,V)= {\cal F}(V,V)$, which is strictly negative by Proposition
\ref{Wirtinger}.

\bigskip

\begin{rema}  The quantity $h''+h+w$ being the radius of  curvature of the curve $\ga,$
 the geometric interpretation of the previous proposition is the following: when
 $h+h''=-w,$ i.e. the radius of curvature vanish, we are at a singularity (vertex) of the curve $\ga,$
 and when $h+h''=w$, the radius of curvature is constant and the corresponding portion of curve
 is an arc of circle of radius $2w.$
\end{rema}

\section{Proof of the main theorem}
Let $(h,w_0(h))$ be a minimizer of ${\cal I}.$ By Proposition
\ref{const}, we know that $|h''+h|$ is constant, and we may assume
without loss of generality that this constant is one. Moreover,
$h''+h$ is characterized, up to multiplication by $-1$ (which does
not change the value of ${\cal F}$) by its set of discontinuities.
The symmetry conditions $(1)$ and $(2)$ imply that there is a
discontinuity at $0,$ and not at $\pi/2.$ Moreover, there must be
at least one discontinuity in the open interval $(0,\pi/2)$:
otherwise $h=A \cos t + B \sin t \pm 1,$ and the boundary
conditions imply $A+(\pm 1)=0$ and $A =0,$ a contradiction. The
remainder of the proof of the main theorem is organized as
follows: we first observe that if $h$ has only one discontinuity
in the interval $(0,\pi/2),$ then the corresponding curve $\ga$ is
the Reuleaux triangle. Then  we prove that if $h$ has at least two
singularities in $(0,\pi/2)$,
 there exists a $h^* \in E$ such
that ${\cal F}(h^*)<{\cal F}(h).$ Therefore there is no minimizer
with at least two singularities and so the only possible one
is the rotated Reuleaux triangle.

\subsection{Case of one interior discontinuity}
 Let $t_1 \in (0,\pi/2)$ be the unique interior discontinuity of $h''+h$. Thus
$$h\left|_{[0,t_{1}]}\right.= A_0 \cos t + B_0 \sin t +1$$
and
$$h\left|_{[t_1,\pi/2]}\right.= A_1 \cos t + B_1 \sin t -1.$$
The conditions $h(0)=0$ and $h'(\pi/2)=0$ imply $A_0=-1$ and $A_1=0.$ On the other hand,
the continuity of $x(t)=(h(t)+w)\cos t- h'(t) \sin t$ at $t_1$ yields
$$ \left\{ \begin{array}{l} A_1= A_0+2 \cos t_1 \\ B_1=B_0, \end{array} \right. $$
so in particular $t_1=\pi/3.$
Changing  the constant $B_0$ (and
thus $B_1$) amounts to making a vertical translation of the curve
$\ga$ and does not affect the geometry of the problem.
Therefore the corresponding curve is unique and  is
nothing but the Reuleaux triangle. From the
computations of Section \ref{s2}, it is easy to compute the volume of the
rotated Reuleaux triangle: we have
$${\cal F}(h_{Reul})=\int_0^{\pi/3} (\cos t -1)dt + \int_{\pi/3}^{\pi/2} \cos t dt = 1-\pi/3.$$
Therefore,
$${\cal I}(B_{Reul})=1+3(1-\frac{\pi}{3})=4-\pi \simeq 0.858407346.$$

\subsection{General case}
Let $0 \leq t_0 <t_1<t_2<  t_3 \leq \pi/2$ and let $h \in E$ such that
\begin{itemize}
\item[-]
$|h''+h|=1$;
\item[-]
  $(t_0,t_1,t_2)$ are three successive discontinuities of $h''+h$;
  \item[-]
 $t_3$ is either the next discontinuity after $t_2,$ or  $t_3=\pi/2.$
\end{itemize}
In particular the case of two interior discontinuities $(t_0,t_3)=(0,\pi/2)$ is covered.
We thus have the following expressions for $h$:
$$h\left|_{[t_0,t_1]}\right. = A_0 \cos t + B_0 \sin t + 1,$$
 $$h\left|_{[t_1,t_2]}\right. = A_1 \cos t + B_1 \sin t - 1,$$
$$h\left|_{[t_2,t_3]}\right. = A_2 \cos t + B_2 \sin t + 1,$$
where $A_i$ and $B_i$ are real constants. As  will become clear
later, the values of the constants $B_0,B_1$ and $B_2$ do not
affect the problem. Next  observe that  the continuity of
$x(t)=(h(t)+w)\cos t- h'(t) \sin t$ at the points $t_1$ and $t_2$
yield the following relations:
$$A_{0} + 2 \cos t_1 = A_1 = A_{2} + 2 \cos t_{2}$$
and thus
$$ \cos t_1 =\frac{A_{1}-A_0}{2} \quad \cos
t_{2}=\frac{A_{1}-A_2}{2}.$$
 From now on we set $x:=-A_0,y:=A_1$ and $z:=-A_2,$ so that $x,y$ and $z$ are three positive constants, and by the assumption
 $t_1<t_2,$ we have $z<x.$ We are going to show that $h$ is not a minimizer of ${\cal F}$,
dividing the proof in three different cases.

\subsubsection{The case $z<y<x$}
  We construct explicitly a map $h^* \in E$ which has one less singularity  than $h$, as follows:
$$ |(h^*)''+h^*|=1,  \forall t \in [0,\pi/2],$$
$$h^*(t)= h(t) , \forall t \in [0,t_0],$$
$$h^*(t)= -h(t) , \forall t \in [t_3,\pi/2],$$
and $(h^*)''+h^*$ has exactly one  discontinuity at $t^* \in
(t_0,t_3).$ Thus
$$h^*\left|_{[t_0,t^*]}\right. = A^* \cos t + B^* \sin t + 1,$$
$$h^*\left|_{[t^*,t_3]}\right. = A^{**} \cos t + B^{**} \sin t - 1.$$
Moreover, as we have $h^*(t_1)=h(t_1)$ and $h^*(t_3)=-h(t_3)$, and
a similar relation for the first derivatives, we deduce that
$A^*=A_0=-x$ and $A^{**}=-A_2=z.$ Finally, the $C^1$ assumption at
$t^*$ implies that

$$ A^{**}=A^* +2 \cos t^* , $$
so that
$$\cos t^* = \frac{A^{**}-A^{*}}{2}=-\frac{A_0+A_2}{2}=\frac{x+z}{2}.$$

\begin{rema}
 If $t_0=0$ and $t_3=\pi/2$, i.e. the case of two singularities,
 one can check that $t^*=\pi/3$, so that $h^*$ corresponds to the Reuleaux triangle.

\end{rema}

Next we compute
$${\cal F}(h^*)-{\cal F}(h)
=\int_0^{\pi/2}((h^*)''+h^*)(h^* \cos t  -(h^*)' \sin t)dt - \int_0^{\pi/2} (h''+h)(h \cos t  -h' \sin t)dt$$
$$ =\int_{t_0}^{t_3}((h^*)''+h^*)(h^* \cos t  -(h^*)' \sin t)dt - \int_{t_0}^{t_3} (h''+h)(h \cos t  -h' \sin t)dt $$
$$
\Big(A^*(t^*-t_{0}) - A^{**}(t_{3}-t^*)\Big)- \Big(
A_{0}(t_{1}-t_{0}) - A_{1}(t_{2} -t_{1})+A_{2}(t_{3}-t_{2})\Big)$$
$$= A_{0}(t^* - t_{1})+A_{1}(t_{2}-t_{1})+A_{2}(t_2-t^*)$$
  $$=(z-x)\arccos \left(\frac{x+z}{2} \right)
+(x-y)\arccos \left(\frac{x+y}{2} \right)
+(y-z)\arccos \left(\frac{y+z}{2} \right).$$

In order to prove that the latter is negative, we first introduce
 the coefficients $a_n$ of the power series  of the function $\arcsin$. It is well known that
$a_n >0, \forall n \geq 1$ and that the radius of convergence of the series is $1.$ Moreover we have
$$\arccos X =\frac{\pi}{2}- \sum_{n=1}^{\infty} a_n X^n.$$
 Next we define the positive map
 $$b_n(a,b,c):=\frac{\left(\frac{a+b}{2}\right)^{n}-\left(\frac{a+c}{2}\right)^{n}}{\frac{b-c}{2}}=
 \left(\sum_{i=0}^{n-1} \left(\frac{a+b}{2} \right)^{i}\left(\frac{a+c}{2} \right)^{n-1-i} \right).$$

Finally we conclude:
$${\cal F}(h^*)-{\cal F}(h)=(z-x)\arccos \left(\frac{x+z}{2}\right)+(x-y)\arccos\left(\frac{x+y}{2}\right)+(y-z)\arccos\left(\frac{y+z}{2}\right)$$
$$=(z-x)\arccos\left(\frac{x+z}{2}\right)+(x-z+z-y)\arccos\left(\frac{x+y}{2}\right)+(y-z)\arccos \left(\frac{y+z}{2}\right)$$
$$=(z-x)\left(\arccos \left(\frac{x+z}{2}\right)-\arccos \left(\frac{x+y}{2}\right)\right)
+(y-z)\left(\arccos \left(\frac{y+z}{2}\right)-\arccos \left(\frac{x+y}{2}\right)\right)$$
$$ =(z-x) \sum_{n=1}^\infty a_n \frac{z-y}{2}b_n(x,z,y)+(y-z)\sum_{n=1}^\infty a_n \frac{z-x}{2}b_n(y,z,x)$$
$$ = \frac{(x-z)(y-z)}{2}\sum_{n=1}^\infty a_n (b_n(x,z,y)-b_n(y,z,x)) < 0,$$
since $b_n(x,z,y)-b_n(y,z,x)$ has the same sign as $y-x.$

\subsubsection{The case $x \leq y$ }

We consider an infinitesimal  variation $h^\eps$ of $h$ such that
$y^\eps=y+\eps$ and we shall prove that $\frac{\pa {\cal F}}{\pa
\eps}(h)<0$. Therefore, by choosing negative $\eps_*$ such that
$|\eps_*|$ is small enough,
 we obtain a map $h^{\eps_*}$ such that ${\cal F}(h^{\eps_*})<{\cal F}(h)$.

From the expressions $\cos t_1=\frac{x+y}{2}$
and $\cos t_2=\frac{y+z}{2}$, we get
$$\cos t_1^\eps = \cos t_1 + \frac{\eps}{2}+o(\eps) \mbox{ and } \cos t_2^\eps = \cos t_2 + \frac{\eps}{2}+o(\eps),$$
which implies
$$ t_1^\eps-t_1= -\frac{1}{2\sin t_1}+o(\eps) \quad \mbox{ and }  \quad  t_2^\eps-t_2=-\frac{1}{2\sin t_2}+o(\eps),$$
while $\cos t_0^\eps=\cos t_0$ and $\cos t_3^\eps = \cos t_3$. In particular, $h^\eps$ and $h$ coincide outside the interval
$(t_0,t_3).$
By a straightforward computation we deduce that
$${\cal F}(h^\eps)-{\cal F}(h)=x(t_1-t_1^\eps) -(y+\eps) (t_2^\eps-t_1^\eps)+y(t_2-t_1)+z(t_2^\eps-t_2)$$
$$=\eps \left( \frac{x}{2\sin t_1} + (t_1-t_2) + y(\frac{1}{2\sin t_2}-\frac{1}{2\sin t_1}) -\frac{z}{2\sin t_2}\right) + o(\eps),$$
hence

$$\frac{\pa {\cal F}}{\pa \eps}(h)=\arccos \left(\frac{x+y}{2}\right)-\arccos \left(\frac{y+z}{2}\right)
   +\frac{x-y}{\sqrt{4-(x+y)^2}}+\frac{y-z}{\sqrt{4-(y+z)^2}}.$$
On the other hand,
$$\arccos \left(\frac{x+y}{2}\right)-\arccos \left(\frac{y+z}{2}\right)
=\int_z^x \frac{\pa \arccos \left(\frac{\xi+y}{2}\right)}{\pa \xi} d\xi < - \frac{x-z}{\sqrt{4-(y+z)^2}},$$
since the map $\xi \mapsto -\frac{1}{\sqrt{4-(y+\xi)^2}} $ is decreasing.
Therefore we have
$$\frac{\partial {\cal{F}}}{\partial\epsilon}(h)<\frac{z-x}{\sqrt{4-(y+z)^2}}+\frac{x-y}{\sqrt{4-(x+y)^2}}+\frac{y-z}{\sqrt{4-(y+z)^2}}$$
$$=\frac{y-x}{\sqrt{4-(y+z)^2}}+\frac{x-y}{\sqrt{4-(x+y)^2}}$$
$$=(y-x)\left(\frac{1}{\sqrt{4-(y+z)^2}}-\frac{1}{\sqrt{4-(x+y)^2}}\right) \leq 0,$$
since $x \leq y $ and using again the fact that $\xi \mapsto -\frac{1}{\sqrt{4-(y+\xi)^2}} $ is decreasing; therefore
$\frac{\partial {\cal{F}}}{\partial\epsilon}(h)<0$.

\subsubsection{The case  $y \leq z$}
As in the previous one, we perform an infinitesimal deformation of
$h^\eps$ of $h$ by setting $y^\eps=y+\eps.$ Here we prove that
$\frac{\partial {\cal{F}}}{\partial\epsilon}(h)>0$ and the
conclusion follows as above. Since the situation is very similar
to the previous case, the details  are left to the Reader.

\bigskip

\noindent
Henri Anciaux, Nikos Georgiou \\
Department of Mathematics and Computing \\
Institute of Technology, Tralee \\
Co. Kerry, Ireland \\
henri.anciaux@staff.ittralee.ie \hspace{2em}
nikos.georgiou@ittralee.ie


\begin{thebibliography}{XXXXX}

\bibitem[Ba]{Ba} T. Bayen, PhD Dissertation, \em Optimisation de formes dans la classe
des corps de largeur constante et des rotors, \em Universit\'e
Pierre et Marie Curie - Paris 6, 2007

\bibitem[BLO]{BLO} T. Bayen, T. Lachand-Robert, E. Oudet, {\em Analytic parametrizations and volume mimization
of three dimensional bodies of constant width}, Arch. Ration.
Mech. Anal. \textbf{186} (2007), no. 2, 225--249

\bibitem[CCG]{CCG} S. Campi, A. Colesanti, P. Gronchi, \em Minimum problems for volumes of
convex bodies, \em in \em  Partial differential equations and applications,
\em Vol. 177 of Lecture Notes in Pure and Appl. Math., Dekker, New York, 1996,  p.
43-55


\bibitem[CG]{CG} G.D. Chakerian and H. Groemer, \em Convex bodies of constant width, \em in \em Convexity and its
applications \em (Ed. P. Gruber and J. Wills) Birkha\"{u}ser, Basel 1983,
49-96



\bibitem[GK]{GK} B. Guilfoyle, W. Klingenberg \em On C$^2$-smooth Surfaces of Constant Width,
\em to appear in Tbilisi Math. Journal

\bibitem[Ha]{Ha} E.  Harrell, II. \em  A direct proof of a theorem of Blaschke and Lebesgue,
\em J. Geom. Anal.
\textbf{12}(2002), no 1, 81-88

\bibitem[Ho]{Ho}
 R. Howard, \em  Convex bodies of constant width and constant brightness, \em  Adv. Math.
\textbf{204}(2006), no 1, 241-261

\end{thebibliography}
\end{document}